\documentclass[11pt]{amsart}

\usepackage{AFBoix2015}

\begin{document}

\author[A.\,F.\,Boix]{Alberto F.\,Boix$^{*}$}
\thanks{$^{*}$Supported by Israel Science Foundation (grant No. 844/14) and Spanish Ministerio de Econom\'ia y Competitividad MTM2016-7881-P}
\address{Department of Mathematics, Ben-Gurion University of the Negev, P.O.B. 653 Beer-Sheva 84105, ISRAEL.}
\email{fernanal@post.bgu.ac.il}

\author[S. Zarzuela]{Santiago Zarzuela$^{**}$}
\thanks{$^{**}$Partially supported by Spanish Ministerio de Econom\'ia y Competitividad MTM2016-7881-P}
\address{Departament de Matem\`atiques i Inform\`atica, Universitat de Barcelona, Gran Via de les Corts Catalanes 585, Barcelona 08007, SPAIN} \email{szarzuela@ub.edu}

\keywords{Frobenius algebras, Cartier algebras, Stanley--Reisner rings, simplicial complexes, free faces.}

\subjclass[2010]{Primary 13A35; Secondary 13F55}

\title[Frobenius and Cartier algebras of Stanley-Reisner rings (II)]{Frobenius and Cartier algebras of Stanley-Reisner rings (II)}

\begin{abstract}
It is known that the Frobenius algebra of the injective hull of the residue field of a complete Stanley--Reisner ring (i.e. a formal power series ring modulo a squarefree monomial ideal) can be only principally generated or infinitely generated as algebra over its degree zero piece, and that this fact can be read off in the corresponding simplicial complex; in the infinite case, we exhibit a 1--1 correspondence between potential new generators appearing on each graded piece and certain pairs of faces of such a simplicial complex, and we use it to provide an alternative proof of the fact that these Frobenius algebras can only be either principally generated or infinitely generated.
\end{abstract}

\maketitle

\dedicatory{To our friend L\^{e} Tu\^{a}n Hoa for his 60th birthday.}

\section{Introduction}\label{introduction}
Let $\Delta$ be a simplicial complex with $n$ vertices, we say
that a pair $(F,G)$ of non--empty, disjoint faces of $\Delta$
is \textbf{free} provided $F\cup G$ is the intersection of all
the facets containing $F.$ Moreover, given two free pairs
$(F,G),(F',G'),$ we say that $(F,G)\leq (F',G')$ if
$F\supseteq F'$ and $G\subseteq G';$ with this partial order, the
set of all the free pairs becomes a partially ordered set. In this
way, a free pair $(F,G)$ is said to be \textbf{maximal} if it is
maximal in the set of free pairs with this order relation.

On the other hand, let $K$ be a field, let $I\subseteq K[x_1,\ldots, x_n]
=R$ be the squarefree monomial ideal attached to $\Delta$ through
the Stanley correspondence, and denote by $I^{[2]}$ the ideal obtaining
after raising to the square all the elements of $I;$ finally, denote
by $J_1$ the smallest ideal of $R$ containing the set
$(I^{[2]}:_R I)\setminus (I^{[2]}+(x_1\cdots x_n)).$

Keeping in mind all the above notations, the main result of this paper
(see Theorem \ref{Josep-Kohji correspondence}) is the below:

\begin{teo}\label{main result: warm up}
There is a 1--1 correspondence between the set of minimal monomial generators
of $J_1$ and the set of maximal free pairs of $\Delta;$ in particular, the number of maximal free pairs of $\Delta$ coincides with the number of minimal monomial generators for $J_1.$
\end{teo}
Our motivation to obtain this result comes from \cite{AlvarezBoixZarzuela2012}, where
the authors focused on the so--called Frobenius and Cartier algebras
of Stanley--Reisner rings; for the convenience of the reader, in what follows
we review some information about these algebras.

Let $A$ be a commutative Noetherian ring and let $M$ be an $A$--module; given any integer $q\geq 1,$ we say that a map $\xymatrix@1{M\ar[r]^-{\phi}& M}$ is \emph{$q$-linear} if, for any $a\in A$ and $m\in M$, $\phi (am)=a^q \phi (m).$ Since the composition of a $q$-linear map with a $q'$-linear map produces a $qq'$-linear map, one can cook up the algebra
\[
F^M:=\bigoplus_{q\geq 1} \End_q (M),
\]
where $\End_q (M)$ denotes the abelian group made up by all the $q$-linear maps on $M;$ the reader will easily note that $F^M$ is an associative, positively $\N$--graded, non--commutative ring, and that its degree $1$ piece is $\End_A (M).$

On the other hand, one says that $\xymatrix@1{M\ar[r]^-{\varphi}& M}$ is \emph{$q^{-1}$-linear} if, for any $a\in A$ and $m\in M$, $\varphi (a^q m)=a \varphi (m).$ Since the composition of $q^{-1}$-linear map with a $q'^{-1}$-linear map produces a $(qq')^{-1}$-linear map one can cook up the algebra
\[
C^M:=\bigoplus_{q\geq 1} \End_{q^{-1}} (M),
\]
where $\End_{q^{-1}} (M)$ denotes the abelian group made up by all the $q^{-1}$-linear maps on $M;$ the reader will easily note that $C^M$ is an associative, positively $\N$--graded, non--commutative ring, and that its degree $1$ piece is $\End_A (M).$

From now on, suppose that $A$ has prime characteristic $p;$ inside $F^M$ and $C^M$ there are respectively distinguished subalgebras, namely
\[
\cF^M:=\bigoplus_{e\geq 0} \End_{p^e} (M),\ \cC^M:=\bigoplus_{e\geq 0} \End_{p^{-e}} (M),
\]
the so--called \emph{algebra of Frobenius (respectively, Cartier) operators} of $M.$ It is known (see either \cite{BlickleBockle2011} or \cite{SharpYoshino2011}) that, if $A=k[\![x_1,\ldots ,x_n]\!]/I$ is a complete local ring having an $F$-finite field $k$ of prime characteristic $p$ as residue field, $R:=k[\![x_1,\ldots ,x_n]\!]$, $E$ denotes the injective hull of $k$ as $R$-module, $E_A$ denotes the injective hull of $k$ as $A$-module, and $(-)^{\vee}:=\Hom_R (-,E)$ denotes the Matlis duality functor, then $\End_{p^e} (E_A)^{\vee}\cong\End_{p^{-e}} (A)$ and $\End_{p^{-e}} (A)^{\vee}\cong\End_{p^e} (E_A).$ From this point of view, under these assumptions, one can think that $\cF^E$ and $\cC^A$ are dual algebras. It is also worth noting that the algebra of Cartier operators, which was introduced by Schwede \cite{Schwede2011} and developed by Blickle \cite{Blickle2013}, has received a lot of attention due to its role in the study of singularities of algebraic varieties in positive characteristic (see \cite{BlickleSchwede2013} and the references therein for additional information).

Hereafter, we will focus on the algebra of Frobenius operators; building upon a counterexample due to Katzman \cite{Katzman2010}, originally motivated by a question raised by Lyubeznik and Smith in \cite{LyubeznikSmith2001}, in \cite{AlvarezBoixZarzuela2012} the authors studied $\cF^{E_A},$ where $A:=K[\![x_1,\ldots ,x_n]\!]/I,$ $K$ is any field of prime characteristic $p,$ $I$ is a squarefree monomial ideal, and $E_A$ denotes the injective hull of $K$ as $A$-module; more precisely, it was proved in \cite[Theorem 3.5 and Remarks 3.1.2]{AlvarezBoixZarzuela2012} that $\cF^{E_A}$ is principally generated as $A$--algebra if and only if $(I^{[2]}: I)=I^{[2]}+(x_1\cdots x_n),$ otherwise it is infinitely generated as $A$--algebra. After \cite{AlvarezBoixZarzuela2012}, the finite or infinite generation of $\cF^{E_A}$ was studied in other situations; for instance, if $A$ is a complete local $F$-finite normal $\Q$--Gorenstein domain of index $m$, then it is known (see \cite[Proposition 4.1]{KatzmanSchwedeSinghZhang2013} and \cite[Theorem 4.6]{EnescuYao2014}) that $\cF^{E_A}$ is finitely generated as $A$--algebra if and only if $m$ is not a multiple of $p$, otherwise it is infinitely generated. In general, keeping in mind that infinite generation of these algebras appear quite often in practice, one can try to understand its growth in this infinite case; with this purpose in mind, Enescu and Yao introduced and studied the so--called \emph{Frobenius complexity}. The interested reader may like to consult \cite{EnescuYao2014} and \cite{EnescuYao2015} for additional information.

One question not answered in \cite{AlvarezBoixZarzuela2012} was whether it is possible to read off the principal (respectively, the infinite) generation of $\cF^{E_A}$ in the simplicial complex $\Delta$ associated to $I$ through the Stanley correspondence; this question was answered in \cite{AlvarezYanagawa2014}, where \`Alvarez Montaner and Yanagawa show that, if $\Delta =\core (\Delta)$ (see Remark \ref{remark about core and cone points}), then $\cF^{E_A}$ is principally generated if and only if $\Delta$ does not have free faces (see \cite[Definition 3.3.1]{Maunder1980} for the definition of free face).

Hereafter, suppose that $\Delta =\core (\Delta),$ and that $\cF^{E_A}$ is infinitely generated as $A$--algebra; on the one hand, by \cite[Theorem 3.5]{AlvarezBoixZarzuela2012}, one knows that $\cF^{E_A}$ is infinitely generated if and only if $(I^{[2]}: I)=I^{[2]}+J_1+(x_1\cdots x_n),$ where $0\neq J_1\not\subseteq I^{[2]}+(x_1\cdots x_n)$ is the smallest monomial ideal containing the set $(I^{[2]}: I)\setminus I^{[2]}+(x_1\cdots x_n).$ On the other hand, by \cite{AlvarezYanagawa2014} $\cF^{E_A}$ is infinitely generated if and only if $\Delta$ has at least one free face. Keeping in mind these two characterizations, one can ask the following:

\begin{quo}\label{asking connection between generators and free faces}
Is there some kind of relation between the number of minimal monomial generators of $J_1$ and the number of free faces of $\Delta$?
\end{quo}

As we will see, such a relation exists but not directly with the free faces of $\Delta$, instead with maximal
free pairs, whereas it is easy to see that free faces are minimal elements of this finite poset.

In fact, the correspondence given in Theorem \ref{main result: warm up} is explicit and one can easily extract from the maximal free pairs the corresponding minimal monomial generators of $J_1$. As application, we use Theorem \ref{main result: warm up} to give an alternative proof of the main technical tool used in \cite[Proof of Theorem 3.5]{AlvarezBoixZarzuela2012} to show that $\cF^{E_A}$ can only be either principally generated or infinitely generated.

Now, we provide a more detailed overview of the contents of this paper for the convenience of the reader. First of all, Section \ref{generators and free faces} contains the main result of this paper (see Theorem \ref{Josep-Kohji correspondence}); namely, a 1--1 correspondence between the non--empty maximal free pairs of a simplicial complex $\Delta$ and the minimal monomial generating set for $J_1.$ Finally, in Section \ref{applications section} we use this correspondence to provide (see Theorem \ref{clarification of the proof of main result in AMBZ12}) an alternative proof of the main technical tool used in \cite[Proof of Theorem 3.5]{AlvarezBoixZarzuela2012} to study the generation of $\cF^{E_A}.$ From here, we immediately conclude (see Theorem \ref{generation of Frobenius algebras and free faces}) that, when $\cF^{E_A}$ is infinitely generated as $A$--algebra, the number of new generators appearing on each graded piece is always upper bounded by a constant, and that this constant is exactly
the number of maximal free pairs of the simplicial complex $\Delta$. We also write down the corresponding dual statement (see Theorem \ref{generation of Cartier algebras and free faces})
for Cartier algebras.

\section{A correspondence between generators and maximal free pairs}\label{generators and free faces}

The purpose of this section is to exhibit a correspondence between minimal monomial generators of $J_1$ and the so--called \textbf{maximal free pairs} of a simplicial complex $\Delta$. Before doing so, we need to state several technical facts; the first one is just an observation made in \cite{AlvarezYanagawa2014} which we use several times in what follows; for this reason, we write it down for the reader's profit (see \cite[Proof of Theorem 4]{AlvarezYanagawa2014} for details). Hereafter, we abbreviate the set $\{1,\ldots ,n\}$ writing just $[n];$ moreover, given a monomial $m\in R$ and an integer $k\geq 0$, set
\[
\supp_k (m):=\{i\in [n]:\ \deg_{x_i} (m)\geq k\}.
\]
The reader will easily note that $\supp_0 (m)=[n]$, $\supp_1 (m)=\supp (m)$, and $\supp_k (m)\subseteq\supp (m)$ for any integer $k\geq 1$; on the other hand, by abuse of notation, hereafter we identify $\supp_k (m)$ with the simplicial complex determined by the elements of $\supp_k (m)$ as vertices. Finally, for any $F\subseteq [n]$ we denote by $\mathbf{x}_F$ the squarefree monomial $\prod_{i\in F}x_i$.

\begin{lm}\label{the free face given by the squares}

Let $m\in R$ be a monomial. Then, $m\notin I^{[2]}$ if and only if $\supp_2 (m)\in\Delta$; moreover, $m\notin I^{[2]}+(x_1\cdots x_n)$ if and only if $\supp_2 (m)\in\Delta$ and $\supp (m)\neq [n]$.
\end{lm}

Next lemma is the first step towards the definition of the so-called \emph{free pairs} (see Definition \ref{generalized free face}).

\begin{lm}\label{towards generalized free faces: step 1}
Let $m\in R$ be a monomial. If $m\in J_1$ and $m \notin I^{[2]}$, then $\supp_2 (m)$ and $[n]\setminus\supp (m)$ are faces of $\Delta$.
\end{lm}

\begin{proof}
By Lemma \ref{the free face given by the squares}, $m\notin I^{[2]}$ if and only if $\supp_2 (m)\in \Delta$ ; so, it only remains to check that $[n]\setminus\supp (m)$ is also a face of $\Delta$.

Indeed, set $G:=[n]\setminus\supp (m)$; assume, to get a contradiction, that $G$ is not a face of $\Delta$. This implies, by the Stanley Correspondence, that $\mathbf{x}_G\in I$ and therefore $m\mathbf{x}_G \in I^{[2]}$, which is equivalent to say that $\supp_2 (m\mathbf{x}_G)\notin\Delta$ (once again by Lemma \ref{the free face given by the squares}). However, since $x_j$ does not divide $m$ for any $j\in G$, one has that $\supp_2 (m)=\supp_2 (m\mathbf{x}_G)\notin\Delta$, a contradiction; the proof is therefore completed.
\end{proof}
It is well--known that, in general, the union of two disjoint faces of a given simplicial complex is not necessarily a face of such a complex; however, in our very specific setting we have the following,  which is the second step towards the definition of free pairs:

\begin{lm}\label{joining faces is also a face}
Let $m\in R$ be a monomial. If $m\in J_1,$ and $m\notin I^{[2]},$ then $\supp_2 (m)\cup [n]\setminus\supp (m)$ is a face of $\Delta$.
\end{lm}

\begin{proof}
Set $H:=\supp_2 (m)\cup [n]\setminus\supp (m),$ and write
\[
m=\left(\prod_{i\in\supp_2 (m)} x_i^{a_i}\right)\prod_{i\notin H}x_i,
\]
where $a_i\geq 2$ are integers. Suppose, to get a contradiction, that $H$ is not a face of $\Delta;$ this implies, by the Stanley Correspondence, that $\mathbf{x}_H \in I$ and therefore
\[
m\mathbf{x}_H=\left(\prod_{i\in\supp_2 (m)} x_i^{a_i+1}\right)\prod_{i\in [n]\setminus \supp_2 (m)}x_i\in I^{[2]}.
\]
This equality shows that, if $M\in I$ is a squarefree generator of $I$ such that $M^2$ divides $m\mathbf{x}_H$, then (by a matter of degrees) $M^2$ will divide $m,$ hence $m\in I^{[2]},$ a contradiction; the proof is therefore completed.
\end{proof}

\begin{rk}\label{connection between free pairs and links}
Let $m$ be a minimal monomial generator of $J_1,$ and set $F:=\supp_2 (m)\neq \emptyset,$ $G:=[n]\setminus\supp (m)\neq \emptyset.$ Lemma \ref{joining faces is also a face} implies, in particular, that $F\in\link_{\Delta} (G)$ and $G\in\link_{\Delta} (F).$
\end{rk}

At this point, recall that a face $F\in\Delta$ is called a \textbf{free face} if $F\cup \{v\}$ is a facet for some $v\notin F$ and $F\cup \{v\}$ is the unique facet containing $F$. We will see now that, for any element $m \in J_1$ such that $m \notin I^{[2]}+(x_1\cdots x_n),$ there exists at least a free face. Although the proof of the below result is verbatim the argument used in \cite[Proof of Theorem 4]{AlvarezYanagawa2014}, we want
to reproduce it here for the sake of completeness.

\begin{lm}\label{towards generalized free faces: step 2}
 Let $m\in R$ be a monomial. If $m\in J_1$ and $m\notin I^{[2]}+(x_1\cdots x_n)$, then there exists $G$
a facet of $\Delta$ such that $\supp_2 (m)\subseteq G,$ and $G$ contains
\textbf{at least} one free face.
\end{lm}

\begin{proof}
Assume, to reach a contradiction, that $\supp_2 (m)$ is contained
in no facet with free faces; our plan is to prove that, in this case, $m\notin J_1$.

First of all, if $\# \supp(m)\leq n-2,$ and $i\notin\supp (m),$ then
$x_i m\notin I^{[2]}+(x_1\cdots x_n),$ and $x_i m\notin (I^{[2]}: I)$ implies
$m\notin (I^{[2]}: I).$ Hence we can replace $m$ by $x_i m$ in this case; repeating
this operation, we may assume that $\# \supp(m)=n-1.$ Let $x_l$ be the unique variable
which does not divide $m.$

Set $F:=\supp_2 (m)\in\Delta,$ we claim that there is a facet $G\supseteq F$ of $\Delta$
with $l\notin G;$ indeed, pick any facet $H\supseteq F$ of $\Delta .$ If $l\notin H,$ then
just take $G=H.$ Otherwise, $H\setminus\{l\}$ \textbf{is contained in a facet $H'$ other
than $H$ by our assumption;} therefore, we can take $G=H'.$ Summing up, replacing $m$ by
\[
\left(\prod_{i\in G\setminus F} x_i\right)m
\]
one can assume that $F=\supp_2 (m)$ is a facet with $l\notin F;$ keeping this in mind, set
\[
M:=x_l\left(\prod_{i\in F} x_i\right).
\]
Notice that $M\in I,$ because $\supp (M)=F\cup\{l\}$
is not a face of $\Delta;$ moreover, since $\supp_2 (mM)
=F\in\Delta,$ Lemma \ref{the free face given by the squares}
ensures that $mM\notin I^{[2]},$ and this finally shows
that $m\notin J_1,$ reaching the contradiction we wanted to show.
\end{proof}

The below result is the third step towards our definition of free pairs. It is also the key technical tool to produce a distinguished family of monomials of $J_1;$ later on (see Theorem \ref{Josep-Kohji correspondence}), we will see that from this family we will extract the minimal
monomial generating set of $J_1.$

\begin{lm}\label{towards generalized free faces: step 3}
Let $m\in R$ be a monomial; suppose that $F:=\supp_2 (m)$
and $G:=[n]\setminus\supp (m)$ are non--empty faces of
$\Delta$ such that $F\cup G\in\Delta,$ and
$F\cup G$ is equal to the intersection of all the facets containing
$F.$ Then, one has that $m\in J_1.$
\end{lm}

\begin{proof}
First of all, since $\supp (m)\neq [n]$ and $\supp_2 (m)\in\Delta,$ Lemma
\ref{the free face given by the squares} ensures
that $m\notin I^{[2]}+(x_1\cdots x_n),$ so it only remains
to show that $m\in (I^{[2]}: I).$

As in the statement, set $F:=\supp_2 (m),$
$G:=[n]\setminus\supp (m),$ and let $M\in I$ be a squarefree monomial, so
$\supp(M)\notin\Delta;$ since $F\cup G\in\Delta$ and
$\supp (M)\notin\Delta,$ for any facet $H\supseteq F$ (and
hence $H\supseteq F\cup G$ by our assumption) there is
$j\in\supp (M)\cap\supp (m)$ such that $j\notin H.$ Now, notice that
$\supp_2 (mM)\supseteq F\cup\{j\};$ thus, if $\supp_2 (mM)\in\Delta$ then
$F\cup\{j\}\in\Delta$ and therefore there is a facet
$H$ of $\Delta$ such that $F\cup\{j\}\subseteq H,$ in particular
$j\in H,$ a contradiction by our choice of $j.$ In this way, one has that
$F\cup\{j\}\notin\Delta,$ hence $\supp_2 (mM)\notin\Delta$ and therefore it follows from Lemma \ref{the free face given by the squares}
that $mM\in I^{[2]}.$ Summing up, $m\in J_1,$ as claimed.
\end{proof}








After all the previous preliminaries, we are now in position to introduce the notion of free pair in the below manner; the reader will easily note that our
inspiration comes from Lemma \ref{towards generalized free faces: step 3}.

\begin{df}\label{generalized free face}
Let $F,G$ be non--empty subsets of $[n].$ On the one hand, we say that $(F,G)$ is a \textit{pair} if $F\cap G=\emptyset,$ and
$F\cup G$ is a face of $\Delta.$ In addition, we say that $(F,G)$ is a \textit{free pair} if $(F,G)$
is a pair, and $F\cup G$ is equal to the intersection of all the facets containing
$F.$
\end{df}

\begin{rk}
Regarding the notion of free pair, we want to single
out the following facts. On the one hand, the reader will easily note (cf.\,Remark \ref{connection between free pairs and links}) that, if $(F,G)$ is a free pair, then $F\in\link_{\Delta} (G)$ and $G\in\link_{\Delta} (F).$ On the other hand, it is also clear that, if $F$ is a free face, and $v\notin F$ is the
vertex such that $F\cup\{v\}$ is the unique facet containing $F$, then the pair $(F,\{v\})$ is free; this fact justifies our choice of terminology.
\end{rk}
Another important observation to keep in mind is that two different free pairs can be contained into the same facet; for instance, pick $I=(xy,xz)$. In this case, the facets of $\Delta$ are $\{1\}$ and $\{2,3\}$; it is also clear that $(\{2\},\{3\})$ and $(\{3\},\{2\})$ are free pairs contained in the same facet.

Motivated by this fact, we introduce the following partial order on the set of free pairs of a given simplicial complex $\Delta$.

\begin{df}\label{poset structure}
Let $\Delta$ be a simplicial complex, and let $\FP (\Delta)$ be the set of all possible free pairs of $\Delta$. Given $(F,G), (F',G')\in \FP(\Delta)$, it is said that $(F,G)\leq (F',G')$ provided $F\supseteq F',$ and $G\subseteq G'.$
\end{df}

Now, we are in position to introduce the following:

\begin{df}\label{maximal gff}
It is said that a free pair on $\Delta$ is \emph{maximal} if it is a maximal element in the poset $\FP (\Delta)$.
\end{df}

\begin{rk}
Notice that if $F$ is a free face, and $v\notin F$ is the
vertex such that $F\cup\{v\}$ is the unique facet containing $F$, then the pair $(F,\{v\})$
is a \textbf{minimal free pair.}
\end{rk}

The reader might ask why we put this specific partial order
on the set of free pairs; the below elementary result is our
response to this question.

\begin{lm}\label{partial order and divisibility}
Given free pairs $(F_1,G_1),(F_2,G_2),$ for each $j=1,2$ set
\[
m_j:=\left(\prod_{i\in F_j}x_i^2\right)\cdot\left(\prod_{i\notin F_j\cup G_j}x_i\right).
\]
Then, $(F_1,G_1)\leq (F_2,G_2)$ if and only if $m_2$ divides $m_1.$
\end{lm}

Next Lemma is the last technical fact we need before showing that
the maximal free pairs of $\Delta$ give rise to a generating
set for $J_1$ (see Lemma \ref{towards generalized free faces: step 4}).

\begin{lm}\label{producing monomials of j1 by divisibility}
If $m\in J_1$ is a monomial, and set $H$ as the biggest subset of $\supp(m)\setminus\supp_2 (m)$ with
the property that $\supp_2 (m)\cup [n]\setminus\supp (m)\cup H\subseteq H'$ for any
facet $H'\supseteq F;$ moreover, suppose that $H \varsubsetneq\supp(m)\setminus\supp_2 (m).$ Then $m/\mathbf{x}_H\in J_1.$
\end{lm}

\begin{proof}
Since $m\in J_1,$ $\mathbf{x}_H m\in J_1.$ Moreover, since $H \varsubsetneq\supp(m)\setminus\supp_2 (m),$ and $\supp_2 (\mathbf{x}_H m)\subseteq
\supp_2(m)\cup H\in\Delta$ (indeed, this is by definition of $H$) and therefore Lemma \ref{the free face given by the squares} ensures that $\mathbf{x}_H m\notin I^{[2]},$ and this implies, by Lemma \ref{joining faces is also a face}, that
$\supp_2 (\mathbf{x}_H m)\cup [n]\setminus\supp (\mathbf{x}_H m)\in\Delta.$ Moreover, since $[n]\setminus\supp (\mathbf{x}_H m)=[n]\setminus\supp (m)$ and
$\supp_2 (\mathbf{x}_H m)=\supp_2 (m)\cup H,$ one has that
\[
(\supp_2 (m))\cup ([n]\setminus\supp (m)\cup H)=(\supp_2 (m)\cup H)\cup [n]\setminus\supp (m)\in\Delta.
\]
This shows that $m/\mathbf{x}_H\in R$ which satisfies by construction, on the one hand, that $\supp_2 (m/\mathbf{x}_H)\cup [n]\setminus\supp (m/\mathbf{x}_H)\in\Delta,$ and, on the
other hand, that
\[
\supp_2 (m/\mathbf{x}_H)\cup [n]\setminus\supp (m/\mathbf{x}_H)=\bigcap_{H'\supseteq\supp_2 (m/\mathbf{x}_H)\text{ facet}}H'.
\]
All these conditions imply, by Lemma \ref{towards generalized free faces: step 3}, that $m/\mathbf{x}_H\in J_1,$ just what we wanted to show.
\end{proof}

We conclude this part by showing that maximal free pairs
of $\Delta$ give rise to a generating set for $J_1;$ more
precisely:

\begin{lm}\label{towards generalized free faces: step 4}
Let $m\in J_1$ be a monomial and assume that $m \notin I^{[2]}+(x_1\cdots x_n)$. Then, there exists $m'\in J_1$ dividing $m$
such that $(\supp_2 (m'),[n]\setminus\supp (m'))$ is a maximal free pair.
\end{lm}

\begin{proof}
Set $F:=\supp_2 (m),$ $G:=[n]\setminus\supp(m),$ and $H$ as the biggest subset
of $\supp (m)\setminus\supp_2 (m)$ with the property that $F\cup G\cup H\subseteq H'$ for any
facet $H'\supseteq F,$ and write
\[
m=\left(\prod_{i\in F}x_i\right)\cdot\mathbf{x}_H\cdot\left(\prod_{i\notin F\cup G\cup H}x_i\right).
\]

Let us show first that $H\varsubsetneq \supp(m)\setminus\supp_2 (m)$ (and so $F\cup G\cup H \neq \emptyset$). Indeed, let $k \in \supp_2(m)$ and $s \in [n]$ such that $x_kx_s \in I$. Then, because $m\notin I^{[2]}$, $s\notin \supp_2(m)$. Assume that $H=\supp(m)\setminus\supp_2 (m)$. This implies that the monomial $\left(\prod_{i\in F}x_i\right)\cdot x_i \notin I$ for any $i\in \supp(m)\setminus\supp_2 (m)$, hence cannot be divided by $x_kx_s$. This implies that $s\neq i$ for any $i\in \supp(m)\setminus\supp_2 (m)$ and so the monomial $mx_kx_s \notin I^{[2]}$. But this is a contradiction with the fact that $m\in J_1$.

Now, Lemma \ref{producing monomials of j1 by divisibility} ensures that
\[
m_1:=\left(\prod_{i\in F}x_i^2\right)\cdot\left(\prod_{i\notin F\cup G\cup H}x_i\right)
\]
is a monomial of $J_1$ which, by construction, satisfies that $(\supp_2 (m_1),[n]\setminus\supp(m_1))$ is a free pair (\textbf{not necessarily maximal}). Now, let $(\supp_2 (m_1),[n]\setminus\supp(m_1))\leq (F',G')$ be a maximal free pair; by Lemma \ref{partial order and divisibility} the monomial
\[
m':=\left(\prod_{i\in F'}x_i^2\right)\cdot\left(\prod_{i\notin F'\cup G'}x_i\right)
\]
divides $m_1$ (\textbf{hence divides $m$}), and belongs to $J_1;$ this is the monomial we are looking for.
\end{proof}

\subsection{Main result}

Now, we are in position to state and prove the main result of this paper, which turns out to be a generalization of the correspondence obtained implicitly during \cite[Proof of Theorem 4]{AlvarezYanagawa2014}.

\begin{teo}\label{Josep-Kohji correspondence}
On the one hand, given a pair $(F,G)$ (see Definition \ref{generalized free face}) set
\[
A(F,G):=\left(\prod_{i\in F} x_i^2\right)\left(\prod_{i\notin F\cup G} x_i\right).
\]
On the other hand, given a monomial $m\in R$ set $Y(m):=(\supp_2 (m),[n]\setminus\supp (m)).$ Then, the set $\{A(F,G):\ (F,G)\text{ maximal free pair of }\Delta\}$
is the minimal monomial generating set for $J_1.$ Moreover, $A$ and $Y$ define $1-1$ correspondences between the set of maximal free pairs of $\Delta$ and the set of minimal monomial generators of $J_1$.
\end{teo}

\begin{proof}
First of all, set $S:=\{A(F,G):\ (F,G)\text{ maximal free pair of }\Delta\}.$ By Lemma \ref{towards generalized free faces: step 4}, $S$ is a monomial
generating set for $J_1,$ so it only remains to show that
is the minimal one; for this, it is enough to check that, if
$(F,G)$ and $(F',G')$ are two different maximal free pairs, then
$m:=A(F,G)$ does not divide $m':=A(F',G').$

Indeed, suppose that $m$ divides $m',$ in particular, for any
$i\in [n],$ $\deg_{x_i} (m)\leq\deg_{x_i} (m');$ this implies
that $F\subseteq F'$ and $G'\subseteq G,$ and so
$(F',G')\leq (F,G),$ hence by the maximality of $(F',G'),$ $F=F'$
and $G=G',$ a contradiction because we start with two \textbf{different} monomials $m$ and $m'.$

Finally, we show that $A$ and $Y$ are bijections being either of them the inverse of the other one; indeed, on the one hand given $(F,G)$ a free pair of $\Delta$, one has that the support of the monomial
\[
m:=\left(\prod_{i\in F} x_i^2\right)\left(\prod_{i\notin F\cup G} x_i\right)
\]
is exactly $[n]\setminus\{G\}$, hence $(Y\circ A) (F,G)=Y(m)=(F,G).$ On the other hand, given now $M$ a minimal monomial generator of $J_1$ not contained in $I^{[2]}+(x_1\cdots x_n)$, one has
\[
(A\circ Y)(M)=A(\supp_2 (M),[n]\setminus\supp (M))=\left(\prod_{i\in\supp_2 (M)} x_i^2\right)\left(\prod_{i\notin\supp_2 (M)\cup ([n]\setminus\supp (M))} x_i\right)=M.
\]
The proof is therefore completed.
\end{proof}



\subsection{Some examples}\label{computing free faces}
One can turn Theorem \ref{Josep-Kohji correspondence} into a naive algorithm which, receiving any simplicial complex $\Delta$ as input, returns all its non--empty maximal free pairs; this method works in the below way.

\begin{enumerate}[(i)]

\item The input is a simplicial complex $\Delta$ with $n$ vertices, and initialize $L$ as the empty list.

\item Compute the corresponding Stanley--Reisner ideal $I\subseteq K[x_1,\ldots, x_n],$ where $K$ is any field.

\item Compute $J:=(I^{[2]}: I)/(I^{[2]}+(x_1\ldots x_n)).$

\item If $J=0,$ then output that $\Delta$ has no free pairs and stop.

\item Otherwise, for each minimal monomial generator $m$ of $J$, add to $L$ the pair $(F,G),$ where
\[
F:=\{i\in [n]:\ \deg_{x_i} (m)=2\}\text{ and }G:=\{i\in [n]:\ \deg_{x_i} (m)=0\}.
\]

\item Output the list $L$.

\end{enumerate}
We have already implemented this algorithm in Macaulay2 (see \cite{M2} and \cite{BoixZarzuelaM2}); next, we show an example where we explain, not only how to use our implementation, but also how to interpret the output.

\begin{ex}\label{example to illustrate the method}
Let $\Delta$ be the simplicial complex of $4$ vertices given by facets $\{1,2\}$ and $\{3,4\};$ it is easy to see that, for each vertex, there is exactly a maximal free pair and these are all. We check this by using our implementation in the below way:
\begin{verbatim}
clearAll;
load "FreePairs.m2";
R=ZZ/2 [x,y,z,w,Degrees=>entries id_(ZZ^4)];
I=ideal(x*z,x*w,y*z,y*w);
L=freePairs(I);
L
{{{4}, {3}}, {{2}, {1}}, {{1}, {2}}, {{3}, {4}}}
\end{verbatim}
Indeed, the Stanley--Reisner ideal attached to $\Delta$ in this case is $I:=(xz,xw,yz,yw)\subseteq K[x,y,z,w],$ and
$(I^{[2]}: I)=I^{[2]}+(x^2zw,xyz^2,xyw^2,y^2zw)+(xyzw),$ this shows that the maximal free pairs of $\Delta$ are $(\{1\},\{2\})$, $(\{3\},\{4\})$, $(\{4\},\{3\})$ and $(\{2\},\{1\}).$ \textbf{Notice that, in this specific example, the maximal free pairs are exactly the free faces of $\Delta.$}
\end{ex}
The reader will easily note that, in this example, every facet gives rise to, at least, one free pair; this is, in general, not always true, as the below example illustrates.

\begin{ex}\label{example to illustrate the method 2}
Let $\Delta$ be the simplicial complex of $4$ vertices given by facets $\{1,4\},$ $\{2,3\}$ and $\{3,4\};$ we use our method to determine all the possible maximal free pairs of $\Delta$ as follows:
\begin{verbatim}
clearAll;
load "FreePairs.m2";
R=ZZ/2 [x,y,z,w,Degrees=>entries id_(ZZ^4)];
I=ideal(x*z,x*y,y*w);
L=freePairs(I);
L
{{{2}, {3}}, {{1}, {4}}}
\end{verbatim}
This shows that, whereas facets $\{1,4\}$ and $\{2,3\}$ gives rise to a free pair, $\{3,4\}$ doesn't give it; indeed, this is because $3$ and $4$ are contained
in more than one facet of $\Delta .$
\end{ex}

We continue by exhibiting an example where the set of free pairs of a simplicial complex is not just made up by maximal free pairs.

\begin{ex}\label{example to illustrate the method 3}
Let $\Delta$ be the simplicial complex of $6$ vertices given by facets $\{1,2,3,4\},$ $\{4,5\}$, $\{4,6\}$ and $\{5,6\};$ we use our method to determine all the possible maximal free pairs of $\Delta$ as follows:
\begin{verbatim}
clearAll;
load "FreePairs.m2";
R=ZZ/2 [x,y,z,w,t,u,Degrees=>entries id_(ZZ^6)];
I=ideal(x*t,x*u,y*t,y*u,z*t,z*u,w*t*u);
L=freePairs(I);
L
{{{3}, {1, 2, 4}}, {{2}, {1, 3, 4}}, {{1}, {2, 3, 4}}}
\end{verbatim}
This shows that this simplicial complex has $(\{3\},\{1,2,4\})$, $(\{2\},\{1,3,4\})$ and $(\{1\},\{2,3,4\})$ as maximal free pairs; the reader will easily note that, for instance, the pair $(\{1,2\},\{3,4\})$ is free, but clearly not maximal. On the other hand, \textbf{although $\{1,2,3\}$ is a free face of $\Delta ,$ the pair
$(\{4\},\{1,2,3\})$ is not even a free pair.}
\end{ex}

We want now to show an example where there is one maximal free
pair $(F,G),$ and $F$ \textbf{is not a singleton set.}

\begin{ex}\label{example to illustrate the method 4}
Let $\Delta$ be the simplicial complex of $5$ vertices given by facets $\{1,2,3\},$ $\{2,4\}$ and $\{3,5\};$ one can check, either by using our method or just
by plotting $\Delta,$ that the maximal free pairs
are $(\{1\},\{2,3\}),$ $(\{2,3\},\{1\}),$ $(\{4\},
\{2\})$ and $(\{5\},\{3\}).$ Notice that
$(F,G)=(\{2,3\},\{1\})$ is a maximal free pair, and
$F$ is not a singleton set.
\end{ex}

In all our above examples, notice that $F\cup G$ was the unique facet
of $\Delta$ containing $F;$ we end up with an example where this
does not happen.

\begin{ex}\label{the Oaxaca example: viva Mexico}
Let $\Delta$ be the simplicial complex of $6$ vertices given by facets $\{1,2,3\},$ $\{1,2,6\}$ and $\{3,4,5\};$ we use our method to determine all the possible maximal free pairs of $\Delta$ as follows:
\begin{verbatim}
clearAll;
load "FreePairs.m2";
R=ZZ/2 [x,y,z,w,a,b,Degrees=>entries id_(ZZ^6)];
I=ideal(x*w,x*a,y*w,y*a,z*b,w*b,a*b);
L=freePairs(I);
L 
{{{6}, {1, 2}}, {{5}, {3, 4}}, {{4}, {3, 5}}, {{2}, {1}}, 
{{1}, {2}}}
\end{verbatim}
Notice that, when $(F,G)$ is either $(\{1\},\{2\})$ or $(\{2\},\{1\}),$ the face
$F\cup G$ turns out to be the intersection of facets $\{1,2,3\}$ and $\{1,2,6\},$ which
are the facets containing $F.$
\end{ex}

\section{Going back to Frobenius and Cartier algebras of Stanley--Reisner rings}\label{applications section}
The goal of this section is to use the correspondence established in
Theorem \ref{Josep-Kohji correspondence} to provide (see Theorem \ref{clarification of the proof of main result in AMBZ12}) an alternative proof of the main technical tool used in \cite[Proof of Theorem 3.5]{AlvarezBoixZarzuela2012}.

\subsection{Preliminary calculations}\label{preliminary warm up}

Let $K$ be a field of prime characteristic $p$, let $R$ be the polynomial ring $K[x_1,\ldots ,x_n]$, and let $I$ be a squarefree monomial ideal with minimal primary decomposition $I=I_1\cap\ldots\cap I_t$ (where each $I_k$ is a face ideal) such that $(x_1,\ldots ,x_n)=I_1+\ldots +I_t$ (this assumption is not a restriction, see \cite[Lemma 3]{AlvarezYanagawa2014} for details). 

\begin{rk}\label{remark about core and cone points}
By the Stanley correspondence, the squarefree monomial ideal $I$ corresponds to a simplicial complex $\Delta$, so the assumption $(x_1,\ldots ,x_n)=I_1+\ldots +I_t$ can be read off in $\Delta$; indeed, let $\mathcal{F} (\Delta)$ be the set of all maximal faces (aka facets) of $\Delta$, let $\CP (\Delta):=\{v\in [n]:\ v\in F \ \forall F\in\mathcal{F} (\Delta)\}$ be the set of cone--points of $\Delta$, and let $\core (\Delta)$ be the restriction of $\Delta$ to the set of vertices not in $\CP (\Delta)$. Keeping in mind the previous notations, one has that $(x_1,\ldots ,x_n)=I_1+\ldots +I_t$ if and only if $\CP (\Delta)=\emptyset$ if and only if $\Delta=\core (\Delta).$
\end{rk}

From now on, given any ideal $J$ of $R$ generated by polynomials $f_1,\ldots ,f_t$, and any integer $a\geq 0$, $J^{[p^a]}$ denotes the ideal generated by $f_1^{p^a},\ldots ,f_t^{p^a}$. It is known (see \cite[Page 168]{AlvarezBoixZarzuela2012} for details) that, for any integer $a\geq 0$, $(I^{[p^a]}: I)=I^{[p^a]}+J_a +(x_1\cdots x_n)^{p^a-1}.$

\begin{disc}
The aim of this discussion is to define the ideal $J_a$ to avoid any misunderstanding; indeed, if $(I^{[p^a]}: I)=I^{[p^a]}+(x_1\cdots x_n)^{p^a-1}$, then $J_a$ is just the zero ideal; so, hereafter, we suppose that $I^{[p^a]}+(x_1\cdots x_n)^{p^a-1}\subsetneq (I^{[p^a]}: I).$ Let $\{m_1,\ldots ,m_t\}$ be a minimal monomial generating set for $(I^{[p^a]}: I);$ we can assume, without loss of generality, that there is $1\leq l\leq t$ such that $m_i\in (I^{[p^a]}: I)\setminus I^{[p^a]}+(x_1\cdots x_n)^{p^a-1}$ for any $1\leq i\leq l.$ In this way, set $J_a :=(m_1,\ldots ,m_l);$ from this definition, it is clear that $J_a$ is a monomial ideal of $R$ contained in $(I^{[p^a]}: I);$ moreover, it is also clear that $J_a$ is the smallest ideal of $R$ containing the set $(I^{[p^a]}: I)\setminus I^{[p^a]}+(x_1\cdots x_n)^{p^a-1}.$
\end{disc}
Now, for any integer $e\geq 1$, following \cite[Page 1]{Katzman2010} set $L_1:=0$ and, for $e\geq 2,$
\begin{equation}\label{truncation of the Frobenius algebra up to degree e}
L_e:=\sum_{\substack{1\leq b_1,\ldots ,b_s\leq e-1\\ b_1+\ldots +b_s=e}} S_{(b_1,\ldots ,b_s)},
\end{equation}
where
\[
S_{(b_1,\ldots ,b_s)}:=\left(I^{[p^{b_1}]}: I\right)\cdot\prod_{k=2}^s \left(I^{[p^{b_1+\ldots +b_k}]}: I^{[p^{b_1+\ldots +b_{k-1}}]}\right).
\]
First of all, observe that, for any $2\leq k\leq s$, flatness of Frobenius \cite[Corollary 2.7 and Corollary 2.8]{Kunz1969} implies that
\[
\left(I^{[p^{b_k}]}: I\right)^{[p^{b_1+\ldots +b_{k-1}}]}=\left(I^{[p^{b_1+\ldots +b_k}]}: I^{[p^{b_1+\ldots +b_{k-1}}]}\right),
\]
so our expression of $L_e$ is the same as the one introduced by Katzman in \cite[Page 1]{Katzman2010}; on the other hand, the reader will also easily note that $\left(I^{[p^{b_1}]}: I\right)=I^{[p^{b_1}]}+J_{b_1} +(x_1\cdots x_n)^{p^{b_1}-1}$ and that, for $2\leq k\leq s$,
\[
\left(I^{[p^{b_1+\ldots +b_k}]}: I^{[p^{b_1+\ldots +b_{k-1}}]}\right)=I^{[p^{b_1+\ldots +b_k}]}+J_{b_k}^{[p^{b_1+\ldots +b_{k-1}}]}+\left(x_1\cdots x_n\right)^{p^{b_1+\ldots +b_k}-p^{b_1+\ldots +b_{k-1}}}.
\]
The first result we want to single out is the following:

\begin{lm}\label{bounding above colon's product}
Fix an integer $e\geq 2$, and let $1\leq b_1,\ldots ,b_s\leq e-1$ be integers such that $b_1+\ldots +b_s=e$. Then, $S_{(b_1,\ldots ,b_s)}\subseteq I^{[p^{b_1+\ldots +b_s}]}+J_{b_1}J_{b_2}^{[p^{b_1}]}\cdots J_{b_s}^{[p^{b_1+\ldots +b_{s-1}}]}+(x_1\cdots x_n).$ In particular,
\[
L_e\subseteq I^{[p^e]}+\sum_{\substack{1\leq b_1,\ldots ,b_s\leq e-1\\ b_1+\ldots +b_s=e}}J_{b_1}J_{b_2}^{[p^{b_1}]}\cdots J_{b_s}^{[p^{b_1+\ldots +b_{s-1}}]}+(x_1\cdots x_n).
\]
\end{lm}

\begin{proof}
By increasing induction on $s\geq 2$. Indeed, for $s=2$ one has
\[
(I^{[p^{b_1}]}:I)\cdot (I^{[p^{b_1+b_2}]}: I^{[p^{b_1}]})\subseteq I^{[p^{b_1+b_2}]}+(IJ_{b_2})^{[p^{b_1}]}+J_{b_1}J_{b_2}^{[p^{b_1}]}+(x_1\cdots x_n).
\]
Moreover, by the very definition of $J_{b_2}$, $IJ_{b_2}\subseteq I^{[p^{b_2}]}$, hence $\left(IJ_{b_2}\right)^{[p^{b_1}]}\subseteq I^{[p^{b_1+b_2}]}$; summing up, combining all these facts it follows that
\[
(I^{[p^{b_1}]}:I)\cdot (I^{[p^{b_1+b_2}]}: I^{[p^{b_1}]})\subseteq I^{[p^{b_1+b_2}]}+J_{b_1}J_{b_2}^{[p^{b_1}]}+(x_1\cdots x_n).
\]
Therefore, the result holds for $s=2$.

Now, suppose that $s\geq 3$ and that, by inductive assumption,
\[
\prod_{k=1}^{s-1} \left(I^{[p^{b_1+\ldots +b_k}]}: I^{[p^{b_1+\ldots +b_{k-1}}]}\right)\subseteq I^{[p^{b_1+\ldots +b_{s-1}}]}+J_{b_1}J_{b_2}^{[p^{b_1}]}\cdots J_{b_{s-1}}^{[p^{b_1+\ldots +b_{s-2}}]}+(x_1\cdots x_n).
\]
In this way, keeping in mind that
\[
\left(I^{[p^{b_1+\ldots +b_s}]}: I^{[p^{b_1+\ldots +b_{s-1}}]}\right)=I^{[p^{b_1+\ldots +b_s}]}+J_{b_s}^{[p^{b_1+\ldots +b_{s-1}}]}+\left(x_1\cdots x_n\right)^{p^{b_1+\ldots +b_s}-p^{b_1+\ldots +b_{s-1}}}
\]
and the inclusion given by the inductive assumption one has
\[
\prod_{k=1}^s \left(I^{[p^{b_1+\ldots +b_k}]}: I^{[p^{b_1+\ldots +b_{k-1}}]}\right)\subseteq I^{[p^{b_1+\ldots +b_s}]}+(IJ_{b_s})^{[p^{b_1+\ldots +b_{s-1}}]}+J_{b_1}J_{b_2}^{[p^{b_1}]}\cdots J_{b_s}^{[p^{b_1+\ldots +b_{s-1}}]}+(x_1\cdots x_n).
\]
Once again, since $IJ_{b_s}\subseteq I^{[p^{b_s}]}$ it follows, by the very definition of $J_{b_s}$, that $(IJ_{b_s})^{[p^{b_1+\ldots+b_{s-1}}]}\subseteq I^{[p^{b_1+\ldots +b_s}]}$. In this way, one finally obtains $S_{(b_1,\ldots ,b_s)}\subseteq I^{[p^{b_1+\ldots +b_s}]}+J_{b_1}J_{b_2}^{[p^{b_1}]}\cdots J_{b_s}^{[p^{b_1+\ldots +b_{s-1}}]}+(x_1\cdots x_n),$ as desired.
\end{proof}
Hereafter, we assume that, given an integer $e\geq 0$, $J_e\neq 0;$ this implies, keeping in mind the description of $J_e$ obtained in \cite[page 168]{AlvarezBoixZarzuela2012}, that any minimal monomial generator of $J_e$ can be written as
\[
\mathbf{x}^{\gamma}=\left(\prod_{i\in\supp_{p^e} (\mathbf{x}^{\gamma})} x_i^{p^e}\right)\left(\prod_{i\in\supp_{p^e-1} (\mathbf{x}^{\gamma})} x_i^{p^e-1}\right),
\]
where this equality just means that the components of $\gamma$ can only be either $0,$ $p^e-1$ or $p^e;$ moreover, these three sets are non-empty and determine uniquely the corresponding minimal monomial generator of $J_e$. The reader will easily note that these sets only depend on the choice of the minimal generator $\mathbf{x}^{\gamma}$, but neither on $e$, nor on $p$ (see \cite[3.1.2]{AlvarezBoixZarzuela2012}).

Next result may be regarded as a non-trivial consequence of the correspondence established in Theorem \ref{Josep-Kohji correspondence}.

\begin{lm}\label{al alba vincero 2}
Let $m,m'\in R$ be two different minimal monomial generators of $(I^{[2]}:I)\setminus I^{[2]}+(x_1\cdots x_n)$ with
$\supp(m)=\supp(m').$ Then, $\supp_2 (m)\not\subseteq\supp_2 (m')$.
\end{lm}

\begin{proof}
Suppose, to obtain a contradiction, that $\supp_2 (m)\subseteq\supp_2 (m'),$ and set
\[
(F,G):=(\supp_2(m),[n]\setminus\supp (m)),\ (F',G'):=(\supp_2(m'),[n]\setminus\supp (m')).
\]
In this way, since by assumption $G=G',$ one has that
$(F',G')\leq (F,G),$  hence by the
maximality of $(F',G'),$ $F=F'$ and $G=G',$ a contradiction because we start
with two \textbf{different} minimal generators $m$ and $m'.$
\end{proof}

Notice that it may happen that
$\supp (m)\subseteq\supp (m')$ (see Example \ref{example to illustrate the method 4}); even more, one
may have that $\supp (m)=\supp (m'),$ as the below example
illustrates.

\begin{ex}\label{even equality of supports}
Let $\Delta$ be the simplicial complex of $4$ vertices given by facets $\{1,3\},$ $\{2,3\}$ and $\{4\};$ one can check that, in this case, $I_{\Delta}=(xy,xw,yw,zw)\subseteq K[x,y,z,w],$ and that
the maximal free pairs of $\Delta$ are $(\{1\},\{3\})$ and
$(\{2\},\{3\}),$ hence $J_1=(x^2yw,y^2xw).$ It is
clear that both $x^2yw$ and $y^2xw$ have the same support.
\end{ex}

Next discussion will play a key role very soon (see Proof of Lemma \ref{Moty argument generalized}). 

\begin{disc}\label{the missing variables determine anything}
Let $e\geq 1$; consider the set of all possible supports for the minimal monomial generators $M\in J_e$ and order it by set inclusion. Then, for any minimal monomial generator $M$ such that $\supp (M)$ is minimal the following holds: for any other minimal monomial generator $M'\neq M$ \textbf{with $\supp (M')\neq\supp (M),$} either $\supp (M)\not\subseteq\supp (M'),$ or $\supp (M)\subsetneq\supp (M').$ Note that these minimal supports are independent of $e$ or $p$. 
\end{disc}

Now, let $\mathbf{x}^{\gamma}$ be any minimal monomial generator of $J_e$ chosen as in Discussion
\ref{the missing variables determine anything} with minimal support; up to permutation of $x_1,\ldots ,x_n$ we can assume, without loss of generality, that there are $1\leq l<r<n$ such that $\mathbf{x}^{\gamma}=(x_1\cdots x_l)^{p^e}(x_{l+1}\cdots x_r)^{p^e-1}.$ The reader will easily note that this is always possible (see paragraph after Proof of Lemma \ref{bounding above colon's product}). Moreover, let $m_1,\ldots, m_k$ be the minimal monomial generators of $J_e,$ and set $G_e\subseteq J_e$ as the ideal generated by the $m_i$'s which satisfies
$\supp (\mathbf{x}^{\gamma})=\supp (m_i);$ finally, given $(i_1,\ldots, i_s)\in\{1,\ldots, k\},$ set
$M:=m_{i_1}m_{i_2}^{p^{b_1}}m_{i_3}^{p^{b_1+b_2}}\cdots m_{i_s}^{p^{b_1+\ldots +b_{s-1}}},$
where $1\leq b_1,\ldots, b_s\leq e-1$ and $b_1+\ldots +b_s=e.$ The reader will
easily note that any monomial generator of $J_{b_1}J_{b_2}^{[p^{b_1}]}\cdots J_{b_s}^{[p^{b_1+\ldots +b_{s-1}}]}$ (see Lemma \ref{bounding above colon's product}) is of this form.

In this way, we are finally ready for proving the below:

\begin{lm}\label{Moty argument generalized}
Preserving all the previous notations and choices, the following assertions hold.

\begin{enumerate}[(i)]

\item $\mathbf{x}^{\gamma}\in J_{b_1}J_{b_2}^{[p^{b_1}]}\cdots J_{b_s}^{[p^{b_1+\ldots +b_{s-1}}]}$ if and only if $\mathbf{x}^{\gamma}\in G_{b_1}G_{b_2}^{[p^{b_1}]}\cdots G_{b_s}^{[p^{b_1+\ldots +b_{s-1}}]}.$

\item If $M$ divides $\mathbf{x}^{\gamma},$ then for any $1\leq t\leq r$ and
for any $2\leq j\leq s,$ $\deg_{x_t} (m_{i_j})=p^{b^j}-1.$

\item $\mathbf{x}^{\gamma}\notin J_{b_1}J_{b_2}^{[p^{b_1}]}\cdots J_{b_s}^{[p^{b_1+\ldots +b_{s-1}}]}$.

\end{enumerate}

\end{lm}

\begin{proof}

Part (i) follows essentially from Discussion \ref{the missing variables determine anything}; indeed, pick $M$ any monomial generator of the product $J_{b_1}J_{b_2}^{[p^{b_1}]}\cdots J_{b_s}^{[p^{b_1+\ldots +b_{s-1}}]}$ and assume that $M$ divides $\mathbf{x}^{\gamma}$. Then $\supp (M)\subseteq\supp (\mathbf{x}^{\gamma})$ and if $M=m_{i_1}m_{i_2}^{p^{b_1}}m_{i_3}^{p^{b_1+b_2}}\cdots m_{i_s}^{p^{b_1+\ldots +b_{s-1}}},$ we have that $\supp (m_{i_j}) \subseteq\supp (\mathbf{x}^{\gamma})$ for any {j}. But keeping in mind Discussion \ref{the missing variables determine anything} this implies that $\supp(\mathbf{x}^{\gamma})=\supp(m_{i_j})$ for any $j$ and so $\mathbf{x}^{\gamma}\in G_{b_1}G_{b_2}^{[p^{b_1}]}\cdots G_{b_s}^{[p^{b_1+\ldots +b_{s-1}}]}.$

Secondly, we want to prove part (ii): indeed, assume, to
reach a contradiction, that for some $1\leq t\leq r$ and some
$2\leq j\leq s,$ $\deg_{x_t} (m_{i_j})=p^{b_j};$ this implies, jointly with
part (i), that
\[
\deg_{x_t}(M)=\sum_{u=1}^s p^{b_1+\ldots +b_{u-1}}\deg_{x_t} (m_{i_u})
\geq p^e +(p^{b_1}+\ldots+p^{b_{j-1}}-1)>p^e,
\]
a contradiction because $\deg_{x_t}(\mathbf{x}^{\gamma})\leq p^e$
and $M$ divides $\mathbf{x}^{\gamma};$ this proves
part (ii).

Finally, we prove part (iii); indeed, suppose that
$\mathbf{x}^{\gamma}\in J_{b_1}J_{b_2}^{[p^{b_1}]}\cdots J_{b_s}^{[p^{b_1+\ldots +b_{s-1}}]},$ this
implies, jointly with part (i), that any $M$ as above
dividing $\mathbf{x}^{\gamma}$ has to be of the form
\[
M=m_{i_1}m_{i_2}^{p^{b_1}}m_{i_3}^{p^{b_1+b_2}}\cdots m_{i_s}^{p^{b_1+\ldots +b_{s-1}}},
\]
and $\supp (\mathbf{x}^{\gamma})=\supp (m_{i_j})$ for any $j.$ Moreover, not all the $m_{i_j}$'s in this product can be equal (otherwise, we
would reach a contradiction again by a matter of degrees), in
particular $m_{i_1}\neq m_{i_j}$ for some $2\leq j\leq s;$ in
this way, Lemma \ref{al alba vincero 2} ensures that there
is $l+1\leq t\leq r$ such that $\deg_{x_t} (m_{i_j})=p^{b_j},$ a
contradiction by part (ii).

Summing up, this shows that $\mathbf{x}^{\gamma}\notin J_{b_1}J_{b_2}^{[p^{b_1}]}\cdots J_{b_s}^{[p^{b_1+\ldots +b_{s-1}}]},$ just what we finally wanted to show.
\end{proof}

\subsection{Third main result}
After all the foregoing calculations, we are definitely in position to prove the third main result of this note, which turns out to be the main tool employed in \cite[Proof of Theorem 3.5]{AlvarezBoixZarzuela2012}; remember that we already defined $L_e$ in \eqref{truncation of the Frobenius algebra up to degree e}.

\begin{teo}\label{clarification of the proof of main result in AMBZ12}
Given any integer $e\geq 0$, and assuming that $J_e\neq 0$, then there is at least a minimal monomial generator of $J_e$ which is not included in $L_e$.
\end{teo}

\begin{proof}
Pick $\mathbf{x}^{\gamma}$ chosen as in Discussion \ref{the missing variables determine anything}. Assume, to get a contradiction, that $\mathbf{x}^{\gamma}\in L_e$; since $L_e$ is a sum of monomial ideals, it follows, combining this fact jointly with Lemma \ref{bounding above colon's product}, that there is some $1\leq b_1,\ldots ,b_s\leq e-1$ with $b_1+\ldots +b_s=e$ such that
\[
\mathbf{x}^{\gamma}\in I^{[p^e]}+J_{b_1}J_{b_2}^{[p^{b_1}]}\cdots J_{b_s}^{[p^{b_1+\ldots +b_{s-1}}]}+(x_1\cdots x_n).
\]
On the one hand, since $\mathbf{x}^{\gamma}\in J_e$, it follows by assumption that $\mathbf{x}^{\gamma}\notin I^{[p^e]}$; on the other hand, since there is at least one variable that does not divide $\mathbf{x}^{\gamma}$, one also has that $\mathbf{x}^{\gamma}\notin (x_1\cdots x_n)$. Finally, Lemma \ref{Moty argument generalized} guarantees that $\mathbf{x}^{\gamma}\notin J_{b_1}J_{b_2}^{[p^{b_1}]}\cdots J_{b_s}^{[p^{b_1+\ldots +b_{s-1}}]}$; in this way, we get a contradiction, which ensures that $\mathbf{x}^{\gamma}\notin L_e$, just what we finally wanted to show.
\end{proof}

\subsection{Consequences about Frobenius and Cartier algebras of Stanley--Reisner rings}

Our final goal is to write down what consequences have Theorem \ref{Josep-Kohji correspondence} and Theorem \ref{clarification of the proof of main result in AMBZ12}
about Frobenius and Cartier algebras of Stanley--Reisner rings.

Indeed, as immediate consequence we obtain the below result about the generation of Frobenius algebras of Stanley--Reisner rings.

\begin{teo}\label{generation of Frobenius algebras and free faces}
Let $K$ be any field of prime characteristic, let $R:=K[\![x_1,\ldots ,x_n]\!]$, let $I=I_{\Delta}\subseteq R$
be a squarefree monomial ideal, let $A:=R/I,$ and let $E_A$ denote the injective hull of $K$ as
$A$-module. If $\cF^{E_A}$ is infinitely generated as $A$-algebra, then the number of
new generators appearing on each graded piece is always less or equal than the number of maximal free pairs of the simplicial complex $\Delta.$
\end{teo}

Because of the duality between Frobenius and Cartier algebras in the $F$-finite case, Theorem \ref{generation of Frobenius algebras and free faces}
has the following immediate consequence; namely:

\begin{teo}\label{generation of Cartier algebras and free faces}
Let $K$ be any \textbf{$\mathbf{F}$--finite field} of prime characteristic, let $R:=K[\![x_1,\ldots ,x_n]\!]$, let $I=I_{\Delta}\subseteq R$
be a squarefree monomial ideal, and let $A:=R/I.$ If $\cC^A$ is infinitely generated as $A$-algebra, then the number of
new generators appearing on each graded piece is always less or equal than the number of maximal free pairs of
the simplicial complex $\Delta.$
\end{teo}

\section*{Acknowledgements}
The authors would like to thank Eran Nevo, Claudiu Raicu and Kevin Tucker for their comments on an earlier draft of this manuscript. Part of this work was done when the first named author visited Northwestern University funded by the CASB fellowship program.

\bibliographystyle{alpha}
\bibliography{AFBoixReferences}

\end{document}